\documentclass[12pt]{article}
\usepackage{latexsym}
\usepackage{amsfonts}
\usepackage{amsmath}
\usepackage{amssymb}
\usepackage{graphicx}
\usepackage{latexsym}
\usepackage{amsfonts}
\usepackage{graphicx}
\usepackage{psfrag}

\input{tcilatex}
\begin{document}

\qquad 

\thispagestyle{empty}

\begin{center}
{\Large \textbf{Positive linear maps and eigenvalue estimates for
nonnegative matrices }}

\vskip0.2inR. Sharma, M. Pal, A. Sharma

Department of Mathematics \& Statistics\\[0pt]
Himachal Pradesh University\\[0pt]
Shimla - 171005,\\[0pt]
India \\[0pt]
email: rajesh.sharma.hpn@nic.in
\end{center}

\vskip1.5in \noindent \textbf{Abstract. \ }We show how positive unital
linear maps can be used to obtain some bounds for the eigenvalues of
nonnegative matrices.



\bigskip

\bigskip

\bigskip

\bigskip

\bigskip

\bigskip

\bigskip

\section{\protect\bigskip Introduction}

\setcounter{equation}{0}Let $\mathbb{M}\left( n\right) $ denote the algebra
of all $n\times n$ complex matrices. An element $A\in M(n)$ with entries $%
a_{ij}$ is nonnegative if $a_{ij}\geq 0$ for all $i$ and $j$ and if $%
a_{ij}>0 $ for all $i$ and $j$, then $A$\ is positive matrix. For any
element $A\in M(n)$ the spectral radius of $A$ is defined as%
\begin{equation}
\rho (A)=\max \left\{ \left\vert \lambda \right\vert :\lambda \in \sigma
(A)\right\}
\end{equation}%
where $\sigma (A)$ is the spectrum of $A$. For a nonnegative matrix $\rho
(A) $ is an eigenvalue of $A$ and it is of basic interest to find the bounds
on $\rho (A)$ in terms of the expressions involving the entries of $A$. For
example, the largest row sum (column sum) of a nonnegative matrix is an
upper bound on $\rho (A),$ we have%
\begin{equation}
\underset{1\leq i\leq n}{\min }\underset{j=1}{\overset{n}{\sum }}a_{ij}\leq
\rho (A)\leq \underset{1\leq i\leq n}{\max }\underset{j=1}{\overset{n}{\sum }%
}a_{ij}
\end{equation}%
and%
\begin{equation}
\underset{1\leq j\leq n}{\min }\underset{i=1}{\overset{n}{\sum }}a_{ij}\leq
\rho (A)\leq \underset{1\leq j\leq n}{\max }\underset{i=1}{\overset{n}{\sum }%
}a_{ij}.
\end{equation}%
For more details see, Horn and Johnson $(2013)$.

Our aim here is to study the inequalities involving eigenvalues of
nonnegative matrices in connection with the positive linear maps. A linear
map $\Phi :\mathbb{M}(n)\rightarrow \mathbb{M}(k)$ is said to be positive if 
$\Phi (A)$ is positive semidefinite $\left( \Phi (A)\geq O\right) $ whenever 
$A\geq O.$ It is called unital if $\Phi (I_{n})=I_{k}.$ In the special case
when $k=1$ the map from $\mathbb{M}(n)$ to $%
\mathbb{C}
$ is called linear functional and it is customary to represent it by the
lower case letter $\varphi .$ A fundamental inequality of Kadison $(1952)$
and its complimentary inequality due to Bhatia and Davis $(2000)$ give the
noncommutative versions of the classical Cauchy-Schwarz inequality $(1821)$
and the Popoviciu inequality $(1935)$. These inequalities involving linear
maps are also important in various other contexts. For instance, these
inequalities are used to derive many interesting bounds pertaining to the
spreads of matrices, see Bhatia and Sharma $(2012,14,16)$. In this paper we
extend this technique further and discuss some inequalities involving
eigenvalues of a nonnegative matrix.

We first derive a lower bound for the sum of the smallest and largest
eigenvalues of a nonnegative symmetric matrix and then use Bhatia-Davis
inequality (2000) to derive a lower bound on $\rho (A)$ in terms of
expressions involving linear maps(Lemma 2.1, Theorem 2.1). This also
provides a lower bound for $\rho (A)$ of a nonnegative matrix (not
necessarily symmetric) and relate it with positive unital linear maps
(Corollary 2.2). This inequality can be used to derive lower bounds for $%
\rho (A)$ in terms of the expressions involving entries of $A$. We
demonstrate some cases here (Corollary 2.3-2.7). We give some examples to
compare our results with the corresponding estimates in the literature
(Example 1.1-1.4)

\section{Main results}

\setcounter{equation}{0}\textbf{Lemma 2.1 }Let $A=\left( a_{ij}\right) \in 
\mathbb{M}(n)$ be a nonnegative symmetric matrix and let its eigenvalues be
arranged as $\lambda _{1}(A)\leq \lambda _{2}(A)\leq ...\leq \lambda
_{n}(A). $ Denote by $b_{jj}$ the $j^{\text{th}}$ smallest diagonal entry of 
$A$. Then%
\begin{equation}
\overset{n-1}{\underset{i=2}{\sum }}\lambda _{i}(A)\leq \overset{n}{\underset%
{i=3}{\sum }}b_{ii}
\end{equation}%
and%
\begin{equation}
\lambda _{1}(A)+\lambda _{n}(A)\geq b_{11}+b_{22}.
\end{equation}

\vskip0.0in\noindent \textbf{Proof} : Let $A=N(A)+D(A)$ where $D(A)=$diag$%
(a_{11},a_{22},...,a_{nn}).$ We rename the diagonal entries $a_{ii}$'$s$ of $%
A$ as $b_{ii}$'s such that $b_{jj}$ is the $j^{\text{th}}$ smallest diagonal
entry of $A$. We replace only first the smallest diagonal entry $b_{11}$ of $%
D(A)$ by the second smallest entry $b_{22}$ and denote the resulting
diagonal matrix by $D(B)$. Then, $D(B)-D(A)$ is a positive semidefinite
matrix as one of its diagonal entry is $b_{22}-b_{11}\geq 0$ and all other
entries are zero. It follows that $B-A$ is also positive semidefinite where $%
B=N(A)+D(B)$.

Let $C=N(A)+D(C)$ where $D\left( C\right) =$ diag($%
b_{22},b_{22},...,b_{22}). $ By Perron-Frobenius theorem, we have%
\begin{equation*}
\lambda _{1}\left( N(A)\right) +\lambda _{n}\left( N(A)\right) \geq 0
\end{equation*}%
and since $\lambda _{i}\left( N(A)\right) =\lambda _{i}(C)-b_{22},$ we get
that%
\begin{equation*}
\lambda _{1}(C)+\lambda _{n}(C)\geq 2b_{22}.
\end{equation*}%
Further, $B-C$ is positive semidefinite therefore $\lambda _{i}(B)\geq
\lambda _{i}(C)$ and hence%
\begin{equation}
\lambda _{1}(B)+\lambda _{n}(B)\geq 2b_{22}.
\end{equation}%
We also have, $\underset{i=1}{\overset{n}{\sum }}\lambda _{i}(B)=$ tr$B=%
\underset{i=1}{\overset{n}{\sum }}b_{ii}.$ Therefore, $(2.3)$ yields%
\begin{equation}
\underset{i=2}{\overset{n-1}{\sum }}\lambda _{i}(B)\leq \underset{i=3}{%
\overset{n}{\sum }}b_{ii}.
\end{equation}%
The inequality $(2.1)$ now follows from $(2.4)$ and the fact that $B-A$ is
positive semidefinite and therefore $\lambda _{i}(A)\leq \lambda _{i}(B)$
for all $i=1,2,...,n$. The inequality $(2.2)$ follows from $(2.1)$ on using
the fact that tr$A=\sum \lambda _{i}(A)=\sum a_{ii}=\sum b_{ii}.$ $\
\blacksquare $\vskip0.2in\noindent \textbf{Theorem 2.1} Let $\Phi :\mathbb{M}%
(n)\rightarrow \mathbb{M}(k)$ be a positive unital linear map. Let $A=\left(
a_{ij}\right) \in \mathbb{M}(n)$ be a nonnegative symmetric matrix and
denote by $b_{jj}$ the $j^{\text{th}}$ smallest of $a_{ii}$'s. Then%
\begin{equation}
\rho (A)\geq \frac{b_{11}+b_{22}}{2}+\left( \Phi (A^{2})-\Phi (A)^{2}\right)
^{\frac{1}{2}}.
\end{equation}%
\textbf{Proof} : Let $\lambda _{\min }\left( A\right) $ and $\lambda _{\max
}\left( A\right) $ respectively denote the smallest and the largest
eigenvalue of $A$. It is clear that $\lambda _{\min }\left( A\right)
I_{k}\leq A\leq \lambda _{\max }\left( A\right) I_{k}$. By Bhatia-Davis
inequality $(2000),$ for any positive unital linear map $\Phi :\mathbb{M}%
(n)\rightarrow \mathbb{M}(k)$, we have%
\begin{equation}
\Phi (A^{2})-\Phi (A)^{2}\leq \frac{\left( \lambda _{\max }\left( A\right)
-\lambda _{\min }\left( A\right) \right) ^{2}}{4}.
\end{equation}%
By Kadison's theorem $(1959)$, $\Phi (A^{2})-\Phi (A)^{2}$ is a positive
semidefinite matrix. For a positive definite matrix $B$ we have $B=U^{\ast }$%
diag$\left( \lambda _{1}(B),\lambda _{2}(B)\text{,}...\text{,}\lambda
_{n}(B)\right) U$ and $B^{\frac{1}{2}}=U^{\ast }$diag$\left( \sqrt{\lambda
_{1}(B)}\text{,}\sqrt{\lambda _{2}(B)}\text{,}...\text{,}\sqrt{\lambda
_{n}(B)}\right) U$ for some unitary $U$ for which $UBU^{\ast }$ is a
diagonal matrix. So, from the inequality $(2.6),$ we have%
\begin{equation}
\frac{\lambda _{\max }\left( A\right) -\lambda _{\min }\left( A\right) }{2}%
\geq \left( \Phi (A^{2})-\Phi (A)^{2}\right) ^{\frac{1}{2}}.
\end{equation}%
Further, $A$ is nonnegative symmetric matrix, Lemma 2.1 therefore ensures
that 
\begin{equation}
\lambda _{\min }\left( A\right) +\lambda _{\max }\left( A\right) \geq
b_{11}+b_{22}.
\end{equation}%
For a nonnegative symmetric matrix, $\lambda _{\max }\left( A\right) =\rho
(A).$ Put $\lambda _{\max }\left( A\right) =\rho (A)$ in $(2.7)$ and $\left(
2.8\right) $ and add the two resulting inequalities, the inequality (2.5)
follows immediately. \ $\blacksquare $

We now extend ($2.5$) for arbitrary nonnegative matrices. Let $\left\vert
A\right\vert =\left( \left\vert a_{ij}\right\vert \right) \in \mathbb{M}(n).$
Then a basic inequality of interest in the present context says that if $B$
and $B-\left\vert A\right\vert $ are nonnegative then%
\begin{equation}
\rho (A)\leq \rho (\left\vert A\right\vert )\leq \rho (B).
\end{equation}%
It also follows from (2.8) that if $A$ and $B-A$ are nonnegative then $\rho
(A)\leq \rho (B).$ See, Horn and Johnson (2013).\vskip0.2in\noindent \textbf{%
Corollary 2.2 }Let $A=\left( a_{ij}\right) \in \mathbb{M}(n)$ be
nonnegative. Let $x_{ij}=\underset{i,j}{\min }\left\{ a_{ij},a_{ji}\right\} $
for all $i$ and $j$ and $X=\left( x_{ij}\right) \in \mathbb{M}(n)$. Then,
for any positive unital linear map $\Phi :\mathbb{M}(n)\rightarrow \mathbb{M}%
(k)$, we have%
\begin{equation}
\rho (A)\geq b_{11}+b_{22}+\left( \Phi (X^{2})-\Phi (X)^{2}\right) ^{\frac{1%
}{2}},
\end{equation}%
where $b_{11}$ and $b_{22}$ respectively denote the first and the second
smallest diagonal entries $A$. \vskip0.2in\noindent \textbf{Proof} : For $%
x_{ij}=\underset{i,j}{\min }\left\{ a_{ij},a_{ji}\right\} ,$ we have $0\leq
x_{kl}\leq a_{kl}$ for all $k,l=1,2,...,n$. It follows that $A$, $X$ and $%
A-X $ are nonnegative matrices. Then , by (2.9), $\rho (A)\geq \rho (X).$
Further, $x_{ji}=\underset{i,j}{\min }\left\{ a_{ji},a_{ij}\right\} =x_{ij}$%
, the matrix $X$ is therefore a symmetric matrix. We now apply Theorem 2.1
to the matrix $X$ and $(2.5)$ gives a lower bound for $\rho (X).$ The
assertions of the corollary now follows from the fact that $\rho (A)\geq
\rho (X).$ \ \ $\blacksquare $\vskip0.2in\noindent \textbf{Corollary 2.3 }%
With notations and conditions as in Corollary 2.2, we have%
\begin{equation}
\rho (A)\geq \Phi (X).
\end{equation}%
\textbf{Proof} : The matrix $X$ is nonnegative and symmetric. Therefore, $%
\rho (X)$ is the largest eigenvalue of $X$. It follows that $\rho (X)I-X$ is
positive semidefinite. This implies that for a positive map $\Phi ,$\ the
matrix $\Phi \left( \rho (X)-X\right) $ is also positive semidefinite.
Further, $\Phi $ is linear and unital therefore we must have$\ \rho (A)\geq
\Phi (X).\ \ \blacksquare $\vskip0.2in\noindent On choosing different
positive unital linear maps in (2.10) and (2.11) we can derive various lower
bounds for $\rho (A).$ We mention a few cases here.

Let $\varphi :\mathbb{M}(n)\rightarrow \mathbb{%
\mathbb{C}
}$ and let $\varphi (X)=\frac{1}{n}\underset{i,j}{\sum }x_{ij}.$ Then $%
\varphi $ is a positive unital linear functional. It then follows from
(2.11) that 
\begin{equation}
\rho (A)\geq \frac{1}{n}\underset{i,j}{\sum }x_{ij},
\end{equation}%
where $x_{ij}=\underset{i,j}{\min }\left\{ a_{ij},a_{ji}\right\} .$

Likewise, for positive unital linear functional $\varphi (X)=\frac{1}{n}%
\underset{i,j}{\sum }x_{ij}$, we have 
\begin{equation}
\rho (A)\geq \frac{x_{ii}+x_{jj}}{2}+x_{ij}.
\end{equation}%
\textbf{Example 1.1} : Marcus and Minc $(1992)$ compare the various bounds
for $\rho (A)$ of the matrix%
\begin{equation*}
A=\left[ 
\begin{array}{ccc}
1 & 1 & 2 \\ 
2 & 1 & 3 \\ 
2 & 3 & 5%
\end{array}%
\right] .
\end{equation*}%
Their best bound is $\rho (A)\geq 5.612.$ The best bound of Wolkowicz and
Styan $(1980)$ for this matrix is $\rho (A)\geq 2.33.$ From $(2.12)$ and $%
(2.13)$, we respectively have $\rho (A)\geq 6.33$ and $\rho (A)\geq 6$. The
actual value of $\rho (A)$ to three decimal places is $7.531$. For the matrix%
\begin{equation*}
B=\left[ 
\begin{array}{ccc}
1 & 1 & 0 \\ 
2 & 1 & 3 \\ 
0 & 3 & 5%
\end{array}%
\right] ,
\end{equation*}%
we respectively have from $(2.12)$ and $(2.13)$, $\rho (B)\geq 5$ and $\rho
(B)\geq 6.$ This also shows that $(2.12)$ and $(2.13)$ are independent.\vskip%
0.2in\noindent Likewise, we can derive several interesting lower bounds for $%
\rho (A)$ on using the inequality $(2.10)$. We mention a few cases here.%
\vskip0.2in\noindent \textbf{Corollary 2.4}. For a nonnegative matrix $%
A=\left( a_{ij}\right) \in \mathbb{M}(n),$ we have%
\begin{equation}
\rho (A)\geq \frac{b_{11}+b_{22}}{2}+\underset{k}{\max }\sqrt{\underset{%
k\neq j}{\underset{k=1}{\overset{n}{\sum }}}x_{jk}^{2}}
\end{equation}%
where $b_{11}$ and $b_{22}$ are respectively the first and the second
smallest diagonal entries of $A$ and $x_{jk}=\min \left\{
a_{jk},a_{kj}\right\} .$\vskip0.2in\noindent \textbf{Proof} : Let $\varphi :%
\mathbb{M}(n)\rightarrow \mathbb{%
\mathbb{C}
}$ be defined as $\varphi (A)=a_{ii}$ for any fixed $i=1,2,...,n$. Then $%
\varphi (A)$ is a positive unital linear functional. So, for $\varphi
(X)=x_{kk}$ we have $\varphi (X^{2})=\underset{j=1}{\overset{n}{\sum }}%
x_{kj}^{2}$ and 
\begin{equation}
\varphi (X^{2})-\varphi (X)^{2}=\underset{j\neq k}{\underset{j=1}{\overset{n}%
{\sum }}}x_{kj}^{2}.
\end{equation}%
We now choose that value of $k$ for which right hand side expression in
(2.15) is maximum and on substituting this value of $\varphi (X^{2})-\varphi
(X)^{2}$ in (2.10) we immediately get (2.14). \ $\blacksquare $\vskip%
0.2in\noindent \textbf{Corollary 2.5 }With notations and conditions as in
Corollary 2.4, we have%
\begin{equation}
\rho (A)\geq \frac{b_{11}+b_{22}}{2}+\underset{i,j}{\max }\left( \underset{%
k\neq i}{\sum }\left\vert x_{ki}\right\vert ^{2}+\underset{k\neq j}{\sum }%
\left\vert x_{kj}\right\vert ^{2}+\frac{\left( a_{ii}-a_{jj}\right) ^{2}}{2}%
\right) .
\end{equation}%
\textbf{Proof} : Let $\varphi (X)=\frac{x_{ii}+x_{jj}}{2},$ for $%
i,j=1,2,...,n$. Then $\varphi (X)$ is a positive unital linear functional.
The inequality (2.16) follows on using arguments similar to those used in
the proof of Corollary 2.4. $\blacksquare $\vskip0.2in\noindent \textbf{%
Corollary 2.6 }With notations and conditions as in Corollary 2.4, we have%
\begin{equation}
\rho (A)\geq \frac{b_{11}+b_{22}}{2}+\sqrt{\frac{\text{tr}X^{2}}{n}-\left( 
\frac{\text{tr}X}{n}\right) ^{2}}
\end{equation}%
where tr$X$ denote the trace of $X.$\vskip0.2in\noindent \textbf{Proof} : We
choose positive unital linear functional $\varphi (X)=\frac{\text{tr}X}{n}$
and use (2.10) and arguments similar to those used in the proof of Corollary
2.4, we immidiately get (2.17). \ $\blacksquare $\vskip0.2in\noindent The
bounds on eigenvalues of a matrix when all its eigenvalues are real as in
case of Hermitian matrices have also been studied in the literature. In
particular, Wolkowicz and Styan (1980) have shown that is eigenvalues of a
matrix are real then%
\begin{equation}
\lambda _{\max }(A)\geq \frac{\text{tr}A}{n}+\frac{1}{\sqrt{n-1}}\left( 
\frac{\text{tr}A^{2}}{n}-\left( \frac{\text{tr}A}{n}\right) ^{2}\right) .
\end{equation}%
In case of nonnegative symmetric matrices when all its diagonal entries are
equal the inequalities (2.17) and (2.18) respectively give%
\begin{equation}
\lambda _{\max }(A)\geq a_{11}+\sqrt{\frac{\text{tr}A^{2}}{n}-\left( \frac{%
\text{tr}A}{n}\right) ^{2}}
\end{equation}%
and%
\begin{equation}
\lambda _{\max }(A)\geq a_{11}+\frac{1}{\sqrt{n-1}}\sqrt{\frac{\text{tr}A^{2}%
}{n}-\left( \frac{\text{tr}A}{n}\right) ^{2}}.
\end{equation}%
It is clear that the inequality $(2.19)$ strengthens the inequality $(2.20)$%
. In general, $(2.17)$ and $(2.18)$ are independent.\vskip0.2in\noindent 
\textbf{Example 1.2}. Let%
\begin{equation*}
A=\left[ 
\begin{array}{cccc}
4 & 0 & 2 & 3 \\ 
0 & 5 & 0 & 1 \\ 
2 & 0 & 6 & 0 \\ 
3 & 1 & 0 & 7%
\end{array}%
\right] ,\text{ \ \ }B=\left[ 
\begin{array}{ccccc}
4 & 1 & 1 & 2 & 2 \\ 
1 & 5 & 1 & 1 & 1 \\ 
1 & 1 & 6 & 1 & 1 \\ 
2 & 1 & 1 & 7 & 1 \\ 
2 & 1 & 1 & 1 & 8%
\end{array}%
\right]
\end{equation*}%
The estimate of Wolkowicz and Styan (1980) gives $\lambda _{4}(A)\geq 7.158$
and $\lambda _{5}(B)\geq 7.449$ while from our estimate $(2.17)$ we have $%
\lambda _{4}(A)\geq 7.8541$ and $\lambda _{5}(B)\geq 7.3983.$ This shows
that for nonnegative symmetric matrices the inequalities $(2.17)$ and $%
(2.18) $ are independent. The actual values of $\lambda _{4}(A)$ and $%
\lambda _{5}(B)$ to three decimal places are $9.376$ and $11.171$ ,
respectively.\vskip0.2in\noindent We now derive a lower bound for the
eigenvalue of a nonnegative symmetric matrix on using the Nagy inequality
(1918) that says that if $x_{1},x_{2},...,x_{n}$ are $n$ real numbers then
for $a\leq x_{i}\leq b$, $i=1,2,...,n,$ we have%
\begin{equation}
\frac{1}{n}\underset{i=1}{\overset{n}{\sum }}x_{i}^{2}-\left( \frac{1}{n}%
\underset{i=1}{\overset{n}{\sum }}x_{i}\right) ^{2}\geq \frac{\left(
b-a\right) ^{2}}{2n}.
\end{equation}%
\textbf{Theorem} : Let $A=\left( a_{ij}\right) \in \mathbb{M}(n)$ be a
nonnegative symmetric matrix and let $\lambda _{\min }(A)$ denote the
smallest eigenvalue of $A$. Then%
\begin{equation}
\lambda _{\min }(A)\geq \underset{i}{\min }a_{ii}-\sqrt{\frac{1}{2}\underset{%
i\neq j}{\overset{n}{\sum }}a_{ij}^{2}}.
\end{equation}%
\textbf{Proof} : Let $Y\in \mathbb{M}(n)$ be such that each of its diagonal
entry equals $\underset{i}{\min }a_{ii}$ and all the offdiagonal entries of $%
Y$ and $A$ are same. Then, $A-Y$ is positive semidefinite and $\lambda
_{\min }(A)\geq \lambda _{\min }(Y).$ On using the inequality $\left(
2.21\right) $, we find that%
\begin{equation}
\frac{\text{tr}Y^{2}}{n}-\left( \frac{\text{tr}Y}{n}\right) ^{2}\geq \frac{%
\left( \lambda _{\max }(Y)-\lambda _{\min }(Y)\right) ^{2}}{2n}.
\end{equation}%
A simple calculation shows that%
\begin{eqnarray}
&&\frac{\left( \lambda _{\max }(Y)-\lambda _{\min }(Y)\right) ^{2}}{2n} 
\notag \\
&=&\frac{\left( \lambda _{\max }(Y)-\frac{\text{tr}Y}{n}\right) ^{2}+\left( 
\frac{\text{tr}Y}{n}-\lambda _{\min }(Y)\right) ^{2}+2\left( \lambda _{\max
}(Y)-\frac{\text{tr}Y}{n}\right) \left( \frac{\text{tr}Y}{n}-\lambda _{\min
}(Y)\right) }{2n}.
\end{eqnarray}%
It is clear that $Y$ is nonnegative symmetric matrix and $\frac{\text{tr}Y}{n%
}=$ $\underset{i}{\min }a_{ii}.$ It follows from the Lemma 2.1 that 
\begin{equation}
\frac{\text{tr}Y}{n}-\lambda _{\min }(A)\leq \lambda _{\max }(A)-\frac{\text{%
tr}Y}{n}.
\end{equation}%
From (2.23)- (2.25), we get that%
\begin{equation}
\frac{\text{tr}Y^{2}}{n}-\left( \frac{\text{tr}Y}{n}\right) ^{2}\geq \frac{2%
}{n}\left( \frac{\text{tr}Y}{n}-\lambda _{\min }(A)\right) ^{2}.
\end{equation}%
The inequality (2.26) implies that%
\begin{equation}
\lambda _{\min }(Y)\geq \frac{\text{tr}Y}{n}-\sqrt{\frac{n}{2}}\sqrt{\frac{%
\text{tr}Y^{2}}{n}-\left( \frac{\text{tr}Y}{n}\right) ^{2}}.
\end{equation}%
Further, tr$Y=\underset{i}{\min }a_{ii}$ and%
\begin{equation*}
\frac{\text{tr}Y^{2}}{n}-\left( \frac{\text{tr}Y}{n}\right) ^{2}=\frac{1}{n}%
\underset{i\neq j}{\sum }a_{ij}^{2}.
\end{equation*}%
The inequality (2.22) then follows from (2.27) $\ \blacksquare $\vskip%
0.2in\noindent Lower bounds for the smallest eigenvalue of a matrix when all
its eigenvalues are real have also been studied in the literature, see
Wolkowicz and Styan (1980). A related inequality says that 
\begin{equation}
\lambda _{\min }(A)\geq \frac{\text{tr}A}{n}-\sqrt{n-1}\sqrt{\frac{\text{tr}%
A^{2}}{n}-\left( \frac{\text{tr}A}{n}\right) ^{2}}.
\end{equation}%
We show by means of the following example that (2.22) and (2.28) are
independent.\vskip0.2in\noindent \textbf{Example 1.3}. Let%
\begin{equation*}
C=\left[ 
\begin{array}{ccc}
2 & 3 & 4 \\ 
3 & 2 & 1 \\ 
4 & 1 & 2%
\end{array}%
\right]
\end{equation*}%
From (2.22) and (2.28) we respectively have $\lambda _{\min }(C)\geq 2-\sqrt{%
13}\cong 1.6056$ and $\lambda _{\max }(C)\geq 2-\sqrt{52}\cong -5.211$ while
from the matrix $A$\ in Example 1.2 we have $\lambda _{\min }(A)\geq 0.258$
and $\lambda _{\min }(A)\geq 0.525.$

\vskip0.2in\noindent \textbf{Corollary 2.7 }Let $A=\left( a_{ij}\right) \in 
\mathbb{M}(n)$ be a nonnegative symmetric matrix. Then the second smallest
eigenvalue of $A$ is less than or equal to the third smallest diagonal entry
of $A$.\vskip0.2in\noindent \textbf{Proof} : Let the eigenvalues of $A$ be
arranged as $\lambda _{1}(A)\leq \lambda _{2}(A)\leq ...\leq \lambda
_{n}(A). $ Then, for any principal submatrix $H$ of $A$\ of order $k$, we
have%
\begin{equation*}
\lambda _{i}(A)\leq \lambda _{i}(H)\leq \lambda _{i+n-k}(A)
\end{equation*}%
for $i=1,2,...,k,$ see Bhatia (1997). Let $B$ be a $3\times 3$ principal
submatrix of $A$ whose largest diagonal entry is $b_{33}.$ Then, on applying
(2.1) to $B$, we get that $\lambda _{2}(A)\leq \lambda _{2}(B)\leq b_{33}.$ $%
\blacksquare $\vskip0.2in\noindent \textbf{Example 1.4}. For the matrix $B$
in Example 1.2 the estimate of Wolcowicz and Styan $\left( 1980\right) $
gives $\lambda _{2}(A)\leq 7.449.$ The third smallest diagonal entry of $B$
is $6.$ So, from our Corollary 2.7 we have a quick and better estimate, $%
\lambda _{2}(A)\leq 6.$

\vskip0.2in\noindent \textbf{Acknowledgements}. The authors are grateful to
Prof. Rajendra Bhatia for the useful discussions and suggestions, and Ashoka
University for a visit in January 2020. The support of the UGC-SAP is also
acknowledged.

\end{document}